\documentclass[a4paper]{amsart}
\usepackage{a4wide}
\usepackage{microtype}
\usepackage{amsmath}
\usepackage{amsthm}
\usepackage{amsfonts}
\usepackage{amssymb}
\usepackage{mathrsfs}
\usepackage{tikz}
\usepackage[hidelinks]{hyperref}
\usepackage{graphicx}
\usepackage{overpic}
\usepackage{todonotes}
\usepackage{enumitem}
\usepackage{array}
\usepackage{pdflscape}
\usepackage{breakcites}

\theoremstyle{plain}
\newtheorem{observation}{Observation}
\newtheorem{conjecture}{Conjecture}
\theoremstyle{definition}
\newtheorem{definition}{Definition}

\theoremstyle{remark}
\newtheorem*{remark}{Remark}

\renewcommand{\Re}{\operatorname{Re}}

\def \RR {\mathbb{R}}
\def \NN {\mathbb{N}}
\def \SS {\mathbb{S}}

\def \CC {\mathbb{C}}
\def \ZZ {\mathbb{Z}}

\def \rPD {\mathrm{rPD}}
\def \rGPD {\mathrm{rGPD}}
\def \mG {\mathrm{mG}}

\def \tG {\mathrm{tG}}
\def \G {\mathrm G}
\def \D {\mathrm D}
\def \tD {\mathrm{tD}}
\def \P {\mathrm P}
\def \H {\mathrm H}

\DeclareMathOperator{\sym}{Sym}

\title{Minimal Twin Surfaces}
\author{Hao Chen}
\address{Mathematical Sciences Research Institute, 17 Gauss Way, Berkeley, 94720 California, USA}
\email{hao.chen.math@gmail.com}
\thanks{A revision of the manuscript was done while the author was in residence
	at the Mathematical Sciences Research Institute in Berkeley, California,
	during the Fall 2017 semester, supported by the National Science Foundation
	under Grant No.\ DMS-1440140.}


\begin{document}

\begin{abstract}
	We report some minimal surfaces that can be seen as copies of a triply
	periodic minimal surface (TPMS) related by reflections in parallel mirrors.
	We call them \emph{minimal twin surfaces} for the resemblance with twin
	crystal.  Brakke's Surface Evolver is employed to construct twinnings of
	various classical TPMS, including Schwarz' Primitive ($\P$) and
	Diamond~($\D$) surfaces, their rhombohedral deformations ($\rPD$), and
	Schoen's Gyroid ($\G$) surface.  Our numerical results provide strong
	evidences for the mathematical existence of $\D$ twins and $\G$ twins, which
	are recently observed in experiment by material scientists.  For $\rPD$
	twins, we develop a good understanding, by noticing examples previously
	constructed by Traizet (2008) and Fujimori and Weber (2009).  Our knowledge
	on $\G$ twins is, by contrast, very limited.  Nevertheless, our experiments
	lead to new cubic polyhedral models for the $\D$ and $\G$ surfaces, inspired
	by which we speculate new TPMS deformations in the framework of Traizet.
\end{abstract}

\maketitle

\section{Introduction}\label{sec:Intro}

In material science, a \emph{crystal twinning} refers to a symmetric
coexistence of two or more crystals related by Euclidean motions.  The simplest
situation, namely the \emph{reflection twin}, consists of two crystals related
by a mirror reflection in the boundary plane.

\begin{figure}[hbt]
	\centering \includegraphics[width=.7\textwidth]{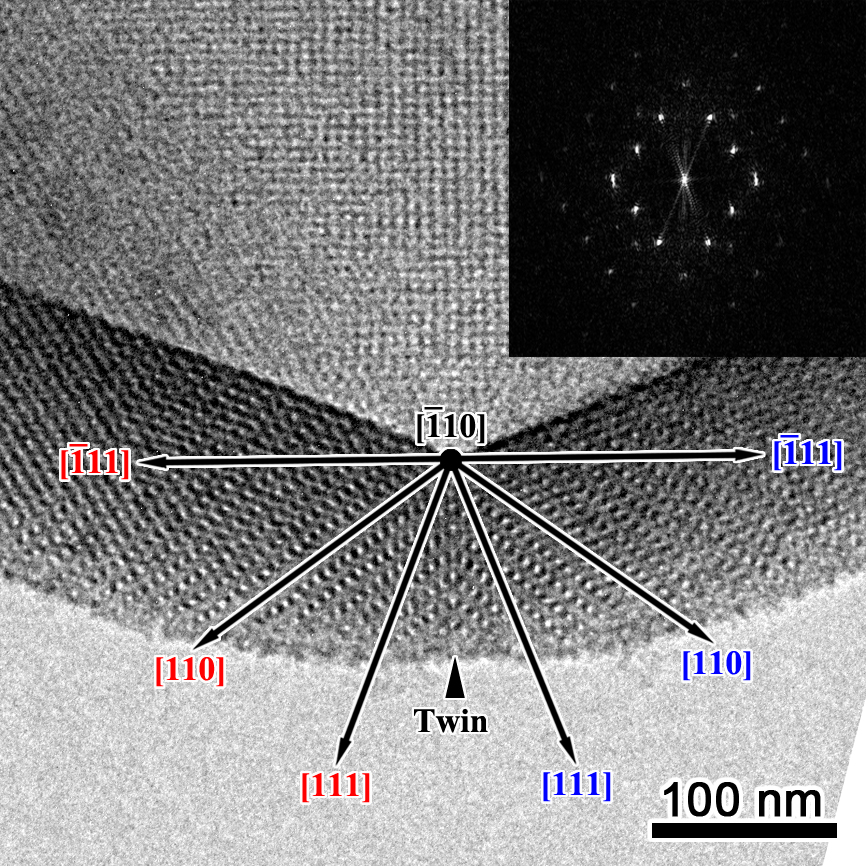}
	\caption{
		\label{fig:han} Twin structure in minimal $\D$ surface observed in
		experiment.  Reuse with permission from \cite{han2011}.  Copyright
		\textcopyright 2011 American Chemical Society.
	}
\end{figure}

\emph{Triply periodic minimal surfaces} (TPMS) are minimal surfaces with the
symmetries of crystals.  They are used to model lyotropic liquid crystals and
many other structures in nature.  Recently, \cite{han2011} synthesized
mesoporous crystal spheres with polyhedral hollows; see also~\cite{lin2017}.  A
crystallographic investigation reveals the structure of Schwarz' $\D$ (diamond)
surface.  Most interestingly, twinning phenomena are observed at the boundaries
of the domains; see Figure~\ref{fig:han}.  We also notice Figure~7.1(b)
in~\cite{hyde1996}, which seems like another evidence, but did not catch the
attention at the time.  Han et al.\ also observed twinnings of Schoen's $\G$
(gyroid) structure\footnotemark, which were earlier discovered
by~\cite{vignolini2012}.

\footnotetext{Paper in preparation.}

However, it is mathematically premature to call the observed structures
``minimal twin surfaces''.  Despite the common belief and various convincing
physics explanations, the energetic base of mesophased systems forming periodic
minimal surfaces is not well understood.  We could not say for sure that the
observed surfaces are minimal.  Moreover, it is \emph{a priori} not known,
mathematically, that a minimal surface with the observed twin structure exists.

\medskip

In this note, we report some minimal surfaces that deserve the name ``minimal
twin surfaces''.  They are similar to \emph{polysynthetic twin crystals},
treating TPMS as crystals.  More specifically, they look like copies of a TPMS
related by mirror reflections in parallel planes.  The minimum lattice distance
$\delta$ between the reflection planes is a parameter of the twin surface.  We
will construct such twinnings for $\rPD$ surfaces and Schoen's $\G$ surface.
$\rPD$ surfaces are rhombohedral deformations of Schwarz' $\P$ and $\D$
surfaces; they are parametrized by a positive real number $t>0$.

We use Brakke's Surface Evolver~\cite{brakke1992}, an efficient gradient
descender, for construction.  A minimal surface is obtained if we manage to
reduce the integral of squared mean curvature down to practically $0$.  Our
main observation is the following:

\begin{observation}
	$\rPD$ twins and $\G$ twins exist for a large set of parameters.
\end{observation}

In particular, $\D$ twins and $\G$ twins exist for very large $\delta$,
providing a strong evidence for the mathematical existence of the twin
structures experimentally observed by~\cite{han2011}.  However, our experiment
for $\rPD$ twins with large parameter $t$ is not conclusive even for small
$\delta$.

\medskip

For the $\rPD$ twins, we notice that extreme examples have been constructed in
previous works.

Examples of small $\delta$ was described by \cite{fujimori2009}.  However, the
period problem was only solved for $\delta=1$ and $\delta=2$.  Following their
work, we observe that, when $\delta=3$, an $\rPD$ twin with sufficiently large
$t$ (near helicoid limit) does not exist.

On the other hand, examples with sufficiently small $t$ (near catenoid limit)
follows from \cite{traizet2008}.  His approach reveals an analogy of $\rPD$ and
$\H$ surfaces with the cubic and hexagonal closed packings.  Further
development of Traizet's technique might lead to rigorous constructions of twin
TPMS with a single reflection plane.

Our understanding of the $\G$ twins is very limited.  However, our numerical
result inspires new cubic polyhedral models for the $\G$ and $\D$ surfaces.
This leads to new interpretations of the $\D$ and $\G$ surfaces in the
framework of \cite{traizet2008}.  Based on numerical evidences, we speculate
new deformations of the $\D$ and $\G$ surfaces.  In particular, we notice a new
tetragonal deformation family with the symmetry of the $\tD$ family.

\medskip

This note is organized as follows.  In Section~\ref{sec:Symmetry} we define
TPMS twinning in analogy with crystal twinning.  In Section~\ref{sec:rPDTwin},
we introduce the $\rPD$ surfaces, define their twins by Weierstrass
representation, and present rigorous examples with small $\delta$
(Section~\ref{ssec:Weber}) and with small $t$~(Section~\ref{ssec:Traizet}), as
well as numerical examples in Surface Evolver (Section~\ref{ssec:rPDevolver}).
The $\G$ twins are numerically constructed in Section~\ref{ssec:Gevolver}.
Then we propose new cubic polyhedral models (Section~\ref{ssec:cubic}) and new
TPMS deformations (Section~\ref{ssec:mG}) for the $\D$ and $\G$ surfaces in the
framework of~\cite{traizet2008}.

\section*{Acknowledgement} 

This manuscript, while self-contained, presents mathematical context and
technical details as a complement and extension to another paper in preparation
on material science.  I would like to thank my collaborators in that project,
especially Chenyu Jin, Lu Han and Nobuhisa Fujita, for their helps on this
note.  I'm very grateful to Karsten Gro{\ss}e-Brauckmann for valuable
suggestions, to Matthias Weber for his comments and Mathematica programs, and
to Han Yu and Martin Traizet for helpful discussions.  Most of the work was
done while the author is visiting St Andrews University, and the author thanks
Louis Theran for his kind invitation.

\section{Crystal twinning and TPMS twinning} \label{sec:Symmetry}

\subsection{Crystals and TPMS}

In physics, a crystal is modeled by a \emph{Bravais lattice}, which can be
defined as a discrete point set $\Lambda \subset \RR^3$ invariant under three
linearly independent translations.  A \emph{lattice plane} of $\Lambda$ is a
plane $H$ such that $H \cap \Lambda$ is a non-empty $2$-dimensional lattice,
i.e.\ invariant under two linearly independent translations.

\medskip 

In this note, the lattice planes of interest are usually placed parallel or
orthogonal to the $z$-axis.  But from time to time, we will use Miller indices
to denote the directions of lattice planes and vectors.  Here is a quick
reference for readers who are not familiar with this notation.

Let $h$, $k$ and $l$ be three coprime integers.  For a lattice spanned by
vectors ${\bf a}, {\bf b}, {\bf c}$
\begin{itemize}
	\item $[hkl]$ denotes the direction of the vector $h{\bf a} + k{\bf b} +
		l{\bf c}$;

	\item $(hkl)$ denotes the lattice planes passing through $n{\bf a}/h$, $n{\bf
		b}/k$ and $n{\bf c}/l$, for some integer $n$;

	\item $\left< hkl \right>$ denotes directions that are equivalent to $[hkl]$
		by symmetry; and

	\item $\{hkl\}$ denotes lattice planes that are equivalent to $(hkl)$ by
		symmetry.
\end{itemize}
Negative numbers in Miller index will be replaced by a number with bar, e.g.\
$\bar 1$ for $-1$.

Alternatively, an $(hkl)$ lattice plane is the set
\[
	H_n = \{ \alpha{\bf a}+\beta{\bf b}+\gamma{\bf c} \in \Lambda \mid
	(\alpha,\beta,\gamma)\in\ZZ^3, \alpha h + \beta k + \gamma l = n \}
\]
for some $n \in \ZZ$.  The \emph{lattice distance} between two $(hkl)$ lattice
planes $H_n$ and $H_m$ is defined as $|m-n|$.  For $H=H_n$, we use $H^+$
(resp.\ $H^-$) to denote the half-space bounded by $H$ containing the parallel
lattice planes $H_m$ with $m>n$ (resp.\ $m<n$).  The normal vector of $H$ is
considered as pointing towards $H^+$.

\medskip

A lattice $\Lambda$ can also be seen as a discrete group acting on $\RR^3$ by
translations; the Bravais lattice is just an orbit of this group.  A
\emph{triply periodic minimal surface} (TPMS) $S \subset \RR^3$ is an
\emph{oriented} minimal surface invariant under the action of a lattice
$\Lambda$.  The quotient $S / \Lambda$ is called the \emph{translational
unit}\footnotemark.  A plane $H$ is a \emph{lattice plane} of $S$ if $H \cap
S$ is invariant under the action of a $2$-dimensional lattice.  Parallel
lattice planes have the same translational symmetries.

\footnotetext{The crystallographic term is ``primitive unit cell''.}

There are several ways to extract a Bravais lattice $\Lambda$ from a TPMS $S$.
For the TPMS that we are interested in (namely $\rPD$ and $\G$ surfaces), we
find it convenient to use flat points as lattice points.  Note that a lattice
plane $H$ for $S$ may not be a lattice plane for $\Lambda$, because $H$ may
contain no flat point at all.  Let ${\bf a}$, ${\bf b}$ and ${\bf c}$ be the
generators of $\Lambda$, then a TPMS lattice plane $H$ can be written as the
set
\[
	H=\{\alpha{\bf a}+\beta{\bf b}+\gamma{\bf c} \mid
	(\alpha,\beta,\gamma)\in\RR^3, \alpha h + \beta k + \gamma l = r \}
\]
for a unique $r \in \RR$.  Note that the coefficients $\alpha$, $\beta$ and
$\gamma$ are in $\RR$ instead of $\ZZ$.  We call the fractional part $\{r\}= r
- \lfloor r \rfloor$ the \emph{offset} of $H$. $H$ contains flat points if and
only if its offset is $0$.  The distance between TPMS lattice planes, as well
as the half-spaces $H^\pm$, are defined in the same way as for a Bravais
lattice. 

\subsection{Twin TPMS and their symmetries}

A lattice plane $H$ of a Bravais lattice $\Lambda$ (resp.\ a TPMS $S$) is said
to be \emph{trivial} if $\Lambda$ (resp.\ $S$) is symmetric under the
reflection in $H$.

The \emph{reflection twin} of a Bravais lattice $\Lambda$ about a non-trivial
lattice plane $H$ consists of $\Lambda \cap H^+$ and its reflective image about
$H$ in the half-space $H^-$.  We call $H$ the \emph{twin boundary}.  Note that
this model is idealized.  In reality, twinning a lattice in this way will
introduce inhomogeneity to the energy, so atoms near the twin boundary will be
slightly pulled away from their original positions.

For a non-trivial $H$ of a TPMS $S$, if we take $S \cap H^+$ and its reflective
image in $\H^-$, the surface obtained is not smooth on the twin boundary, hence
definitely not a minimal surface.

\medskip

Instead of naively imitating the construction for crystal twinning, we will
define TPMS twinning by specifying the symmetry.  Like Bravais lattices, TPMS
often admit extra symmetries other than the translations.  The \emph{symmetry
group} or \emph{space group} of a TPMS $S$, denoted by $\sym(S)$, is the group
of Euclidean symmetries of $S$.  For a lattice plane $H$ of $S$, $\sym(H)$ is
the (setwise) stabilizer of $H \cap S$ in $\sym(S)$.  In other words, $\sym(H)$
is the subgroup of $\sym(S)$ leaving $H \cap S$ (setwise) invariant.  For
example, the following elements of $\sym(S)$ belongs to $\sym(H)$:
\begin{itemize}
	\item the translations in directions parallel to $H$,

	\item the rotations about axis orthogonal to $H$, and

	\item the reflections in a plane orthogonal to $H$.
\end{itemize}
Parallel lattice planes have the same symmetry.  Crystal twinning breaks
symmetries of the lattice, but preserves those symmetries that leave the twin
boundary invariant.  We expect the same for the TPMS twinning.

We say that a surface $S'$ is \emph{asymptotic to} $S$ in a half-space $H^+$ if
there is a sequence of translations $(T_k)_{k \in \NN}$ in the normal direction
of $H$ such that $(S' \cap T_k(H^+))_{k \in \NN}$ converges to $T_k(S \cap
H^+)$.

We are now ready to define the TPMS twinning:

\begin{definition}
	Let $S \subset \RR^3$ be a triply periodic minimal surface, and $H$ be a
	non-trivial lattice plane of $S$.  A \emph{reflection twin} of $S$ about $H$
	is a doubly periodic minimal surface $\Sigma$ that is
	\begin{itemize}
		\item invariant under the action of $\sym(H)$;

		\item invariant under the reflection in $H$;

		\item asymptotic to $S$ in $H^+$.
	\end{itemize}
\end{definition}

\begin{remark}
	The anonymous referee points out that the reflection twins are examples of
	the following general problem: Given two minimal surfaces $S^+$ and $S^-$ and
	a hyperplane $H \subset \RR^3$, find a minimal surface $\Sigma$ that is
	asymptotic to $S^+$ in $H^+$ and asymptotic to $S^-$ in $H^-$.  Earlier
	examples include the Scherk surface, for which $S^+$ is parallel horizontal
	planes, and $S^-$ is parallel vertical planes.
\end{remark}

A reflection twin has a single symmetry plane $H$, hence not periodic in the
direction orthogonal to $H$.  Available tools for minimal surfaces of infinite
topology are very limited.  The author is aware of~\cite{morabito2012}
and~\cite{traizet2013} which, if combined with the arguments
in~\cite{traizet2008}, could lead to a rigorous treatment of reflection twin.
This is the topic of a future project.

The current manuscript focuses on numerical experiments.  Then the aperiodicity
becomes very inconvenient as computers do not really understand infinity.
Therefore, we will work with \emph{polysynthetic twins} with parallel
reflection planes, as an approximation of single reflection twins.

\begin{definition}
	Let $S \subset \RR^3$ be a triply periodic minimal surface, and $H$ be a
	non-trivial lattice plane of $S$.  A \emph{polysynthetic twin} of $S$ (or an
	\emph{$S$ twin}) about $H$ is a triply periodic minimal surface
	$\Sigma_\delta$, for some integer parameter $\delta>0$, such that
	\begin{itemize}
		\item $\Sigma_\delta$ is invariant under the action of $\sym(H)$;

		\item $\Sigma_\delta$ is invariant under the reflection in $H$ and parallel planes;

		\item the minimum lattice distance between reflection planes is $\delta$;

		\item $\Sigma_\delta$ between adjacent reflection planes is ``similar'' to $S$.
	\end{itemize}
\end{definition}

Here the word ``similar'' is not well-defined in general.  Since our goal is to
approximate the reflection twin, we would expect the following asymptotic
behavior for the sequence $(S_\delta)_{\delta \in \NN}$: There is a sequence of
translations $(T_\delta)_{\delta \in \NN}$ in the normal direction of $H$ such
that $(S_\delta \cap T_\delta(H^+))_{d \in \NN}$ converges to $T_\delta(S \cap
H^+)$.  In our construction, this expectation is reflected either in the
proposed Weierstrass representation (see Definition~\ref{def:rPDTwin}), or by
the design of the Surface Evolver experiment (see
Sections~\ref{ssec:rPDevolver} and \ref{ssec:Gevolver}).

In the following, the word ``polysynthetic'' will be omitted unless otherwise
stated.

\section{\texorpdfstring{$\rPD$}{rPD} twin surfaces} \label{sec:rPDTwin}

\subsection{The rPD surfaces} \label{ssec:rPD}

The $\rPD$-family is a $1$-parameter family of TPMS with rhombohedral
symmetries, generalizing Schwarz' $\D$ and $\P$ surfaces.  A Weierstrass
representation for the translational unit of an $\rPD$ surface is given
in~\cite{fogden1993} and \cite{fogden1999}
\[
	S_t \colon \omega \mapsto \Re\int^\omega (1-z^2, i(1+z^2), 2z) R_t(z)\,dz,
\]
where $\omega \in \CC$,
\[
	R_t(z)=[z(z^3-t^3)(z^3+t^{-3})]^{-1/2},
\]
and $t>0$ is the parameter.  The Schwarz' $\D$ and $\P$ surfaces are restored
with $t=\sqrt{1/2}$ and $t=\sqrt{2}$, respectively.  The $\rPD$-family is
self-conjugate; the conjugate surface of $S_t$ is $S_{1/t}$.

The $\rPD$ surfaces are known to \cite{schwarz1890} and rediscovered by
\cite{schoen1970} and \cite{karcher1989}.  They belong to Meeks'
$5$-dimensional family~\cite{meeks1990}.  In the Weierstrass representation of
a Meeks' surface, the Gauss map represents the translational unit as a
two-sheeted cover of $\SS^2$ with four antipodal pairs of simple branch points
(eight in total).  The branch points correspond to the flat points on the
surface.  For instance, after a stereographic projection onto the complex plane
$\CC$, the branch points of $S_t$ are the roots of $R_t(z)^{-2}$.

The Weierstrass data reflects the symmetries of the $\rPD$ surfaces.  More
specifically, the surface has vertical symmetry planes (parallel to the
$z$-axis).  The group generated by these reflections is the Euclidean triangle
group with parameters $(3,3,3)$.  \cite{weyhaupt2006} proved that an embedded
TPMS of genus $3$ with these symmetries must be an $\rPD$ surface or an $\H$
surface.

\medskip

We recommend the following way to visualize $\rPD$ surfaces.  Consider two
equiangular triangles that intersect the $z$-axis perpendicularly at their
centers, whose projections on the $xy$ plane differ by a rotation of $\pi/3$.
Our building block, which we call \emph{catenoid unit}, is the ``catenoid''
spanned by the two triangles.  The whole surface is obtained by the order-$2$
rotations about the edges of the triangles.  The $z$-axis is actually in the
$[111]$-direction of the rhombohedral lattice.

\begin{remark}
	If the projections of the two triangles coincide on the $xy$ plane, then the
	same construction yields an $\H$ surface.  Catenoid units of $\H$ surfaces
	will appear on the twin boundaries of $\rPD$ twins.
\end{remark}

\begin{figure}[hbt]
	\centering
	\includegraphics[width=.7\textwidth]{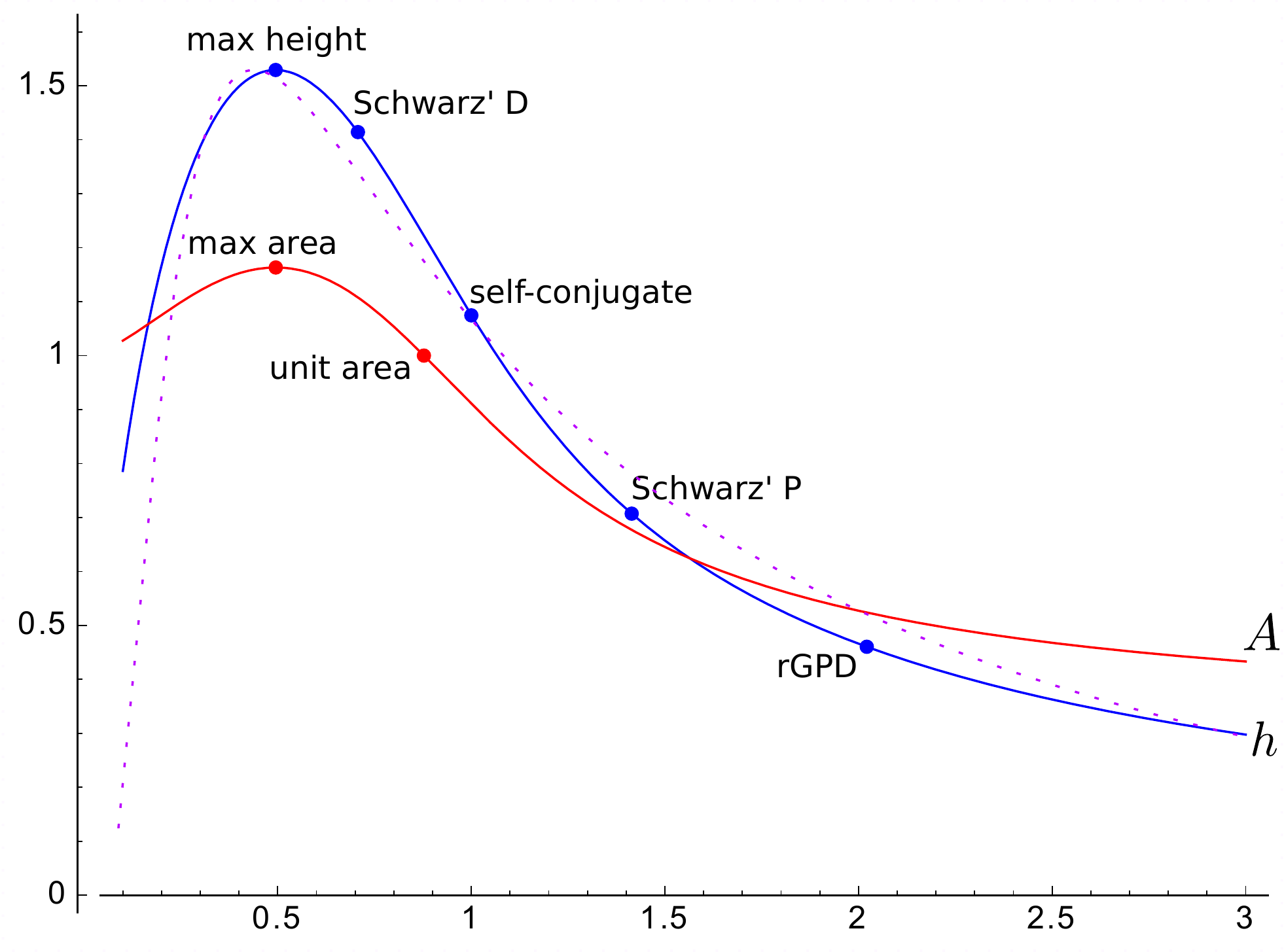}
	\caption{
		Plot of the height $h$ (solid blue) of a catenoid unit assuming unit
		inradii for the bounding triangles, and the area $A$ (solid red) assuming
		unit total area for the bounding triangles, against the parameter $t$ on
		the horizontal axis.  The dotted line is the reparameterization of $h$ with
		the parameter $\tau$ in place of $t$ on the horizontal axis.
		\label{fig:plots} 
	}
\end{figure}

The $1$-parameter family $\rPD$ is obtained by ``stretching'' the two triangles
along the $z$-axis.  Let $h(t)$ denote the height of the catenoid unit
(vertical distance between the triangles) assuming unit inradii for the
triangles, and $A(t)$ be the area of the catenoid assuming unit total area for
two triangles.  Using the formulae given in Appendix A of~\cite{fogden1999}
(where the parameter $r_0$ correspond to our $1/t$), we plot $h$ and $A$
against $t$ in Figure~\ref{fig:plots}.  The calculations are done in
Sage~\cite{sagemath}.

The height $h$ attains its maximum $h_{\max}=1.529295\cdots$ at
$t_0=0.494722\cdots$ and converges to $0$ in both limits $t \to \infty$ and $t
\to 0$.  With a distance bigger than $h_{\max}$, the triangles span no
catenoid; in this case the only minimal surface is the relative interiors of
the triangles.  Any $h<h_{\max}$ corresponds to two different $\rPD$ surfaces.

The area $A$ attains its maximum $A_{\max}=1.163261\cdots$ at
$t=0.494893\cdots$ and converges to $1$ in the limit $t \to 0$, to $1/3$ in the
limit $t \to \infty$.  The area exceeds $1$ as long as $t<t_1=0.877598\cdots$.
For $t$ between $t_0$ and $t_1$, the catenoid is only a local minimizer for the
area functional; the relative interiors of the triangles has smaller area.

Other interesting points are: $t=\sqrt{1/2}$ for the $\D$ surface; $t=\sqrt{2}$
for the $\P$ surface; $t=1$ is a self-conjugate surface; and $t=1/t_0$ is the
$\rGPD$ surface according to~\cite{fogden1999}.

\subsection{A Weierstrass representation of rPD twins} \label{ssec:Weierstrass}

\begin{figure}[hbt]
	\centering
	\includegraphics[width=.6\textwidth]{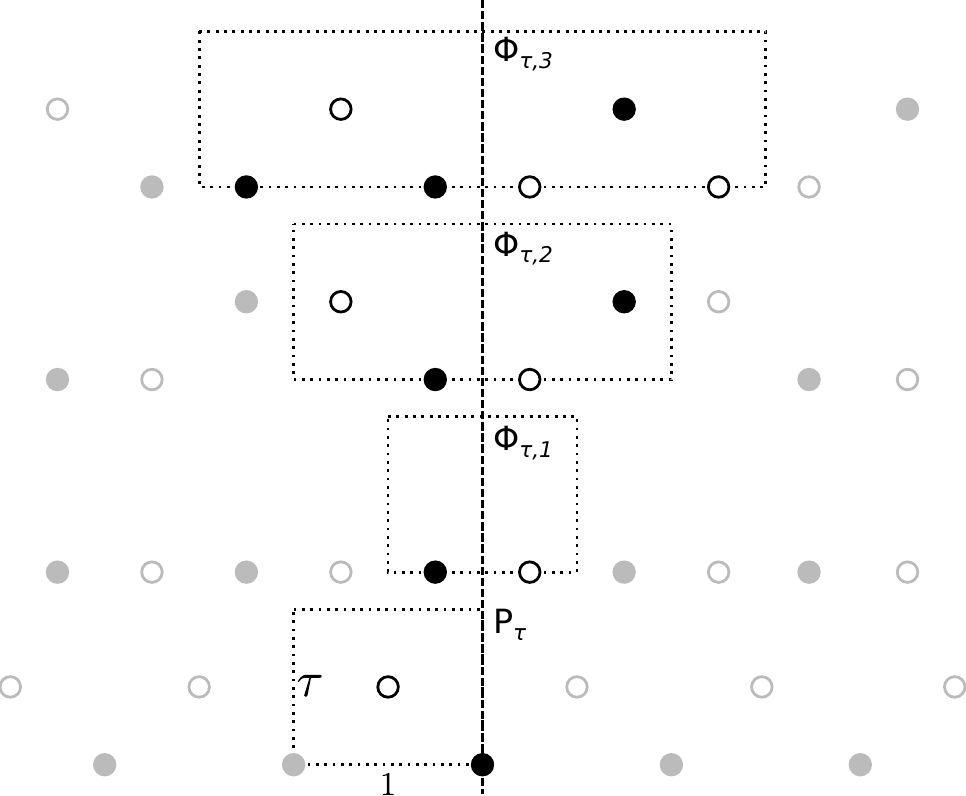}
	\caption{
		\label{fig:Weber} Zeros (filled circles) and poles (empty circles) of the
		Gauss maps $G_{\tau}$ and $G_{\tau,\delta}$, $\delta=1, 2, 3$.  The dotted
		rectangles indicate the tori on which the Gauss maps are defined.
	}
\end{figure}

\cite{fujimori2009} developed an approach to construct surfaces with vertical
symmetry planes.  Instead of a cover of $\SS^2$, the Gauss map descends to the
quotient torus $\CC / \left< 1, \tau i \right>$, where $\tau \in \RR_+$ is an
adjustable parameter.  For example, the Weierstrass representation of an $\rPD$
surface is given by~\cite{weyhaupt2006}\cite{weyhaupt2008}
\[
	S_\tau \colon \omega \mapsto \Re\int^\omega(\frac{1}{2}(G_\tau^{-1}-G_\tau),
	\frac{i}{2}(G_\tau^{-1}+G_\tau), 1)\,dz,
\]
where the Gauss map
\[
	G_\tau(z) = \Big( \rho \frac{\vartheta(z; \tau)}{\vartheta(z-1/2-\tau/2;
	\tau)} \Big) ^{2/3}
\]
is defined on the torus $\CC / \left< 1, \tau i \right>$.  Here, $\rho$ is the
Lopez--Ros factor, and
\[
	\vartheta(z;\tau)=\sum_{k=-\infty}^\infty
	e^{ - \pi (k+\frac{1}{2})^2 \tau + 2 \pi i (k + \frac{1}{2})(z-\frac{1}{2})}
\]
is one of the Jacobi $\vartheta$-functions.  $\vartheta(z;\tau)$ has simple
zeros at the lattice points spanned by $1$ and $\tau i$.  At the bottom of
Figure~\ref{fig:Weber}, we show the zeros and poles of $G_\tau$.  They
correspond to points of $S_\tau$ with vertical normal vectors; these are
exactly the flat points of $S_\tau$, and lie on the intersections of the
vertical symmetry planes.

Note that the parameter $\tau$ is different from $t$ but equivalent.  For
comparison, the height of the catenoid unit is plotted against $\tau$ by a
dotted curve in Figure~\ref{fig:plots}.  The plot is generated by
Sage~\cite{sagemath}, but the calculations are done using mpmath~\cite{mpmath}.
The Schwarz' $\P$ surface is restored with $\tau =
1.563401\cdots$~\cite{weyhaupt2006}.  The conjugate surface of $S_\tau$ is
$S_{1/\tau}$.  We will take the liberty to switch between the two parameters
whenever convenient.

\medskip

Apart from the vertical reflection planes, horizontal reflection planes are
assumed in~\cite{fujimori2009} in order to take advantage of the symmetry.  Up
to a translation, we may assume that one horizontal symmetry plane correspond
to the imaginary axis of $\CC$.

Following the argument in~\cite{fujimori2009}, we propose the following
definition of $\rPD$ twins through Weierstrass representation.

\begin{definition}\label{def:rPDTwin}
	The \emph{$\rPD$ twin surface} $\Sigma_{\tau,\delta}$ is given by the Weierstrass
	representation
	\[
		\Sigma_{\tau,\delta} \colon \omega \mapsto \Re\int^\omega\Big(\frac{1}{2}(G_{\tau,\delta}^{-1}-G_{\tau,\delta}), \frac{i}{2}(G_{\tau,\delta}^{-1}+G_{\tau,\delta}), 1\Big)\,dz,
	\]
	where the Gauss map
	\[
		G_{\tau,\delta}(z) = \rho \prod_{k \text{ odd}} \Big(\frac{\vartheta((z-p_k)/\delta;
		\tau/\delta)}{\vartheta((z+p_k)/\delta; \tau/\delta)}\Big)^{2/3} \prod_{k \text{ even}}
		\Big(\frac{\vartheta((z-p_k-\tau/2)/\delta; \tau/\delta)}{\vartheta((z+p_k-\tau/2)/\delta;
		\tau/\delta)}\Big)^{-2/3},
	\]
	and $(p_k)_{1 \le k \le \delta}$ are real numbers such that $0 < p_1 < p_2 < \dots
	< p_\delta < \delta/2$ and $p_k + p_{\delta+1-k} = \delta/2$.
\end{definition}

Note that we scale the $\vartheta$-function by $1/\delta$ so that, within the
stripe $0 < \Re(z) < \delta/2$, the zeros and poles of $G_{\tau,\delta}$ are
similarly arranged as for $G_\tau$; see Figure~\ref{fig:Weber}.  The Lopez--Ros
factor $\rho$ is fixed to $1$ so that $|G_{\tau,\delta}(z)|=1$ for all $z \in
i\RR$; see Appendix A of~\cite{weyhaupt2006}.

\subsection{Example and non-examples of small \texorpdfstring{$\delta$}{\it d}} \label{ssec:Weber}

The fact that $\Sigma_{\tau,\delta}$ satisfies the expected symmetries follows
from~\cite{fujimori2009}.  But to prove their existence, we need to solve the
period problem.  In our case, the period problem for $\Sigma_{\tau,\delta}$ asks to
find $p_k$ such that
\[
	\Re \int^{p_k}(\frac{1}{2}(G_{\tau,\delta}^{-1}-G_{\tau,\delta}),
	\frac{i}{2}(G_{\tau,\delta}^{-1}+G_{\tau,\delta}))\,dz
\]
are vertices of an equiangular triangle.  In~\cite{fujimori2009}, the period
problem is solved for small $\delta$'s, as an answer is immediate by symmetry.  When
$\delta=1$ we have $p_1=1/4$; this is an $\H$ surface.  When $\delta=2$, we have $p_1=1/4$
and $p_2=3/4$; this is Karcher's T-WP surface.

For $\delta=3$, the period problem is $1$-dimensional.  Examples of $\Sigma_{t,3}$
are computed numerically.  The Mathematica program for this purpose is kindly
provided to us by Matthias Weber.  The period problem asks to solve the
following equations:
\begin{align}\label{eqn:period}
	\begin{split}
		\Re \int_0^{3/2+\tau/2} \frac{i}{2}(G_{\tau,3}^{-1}+G_{\tau,3})\,dz =& 0\\
		\cos\frac{\pi}{6} \Re \int_{\tau/2}^{3/2} \frac{1}{2}(G_{\tau,3}^{-1}-G_{\tau,3})\,dz
		+ \sin\frac{\pi}{6} \Re \int_{\tau/2}^{3/2} \frac{i}{2}(G_{\tau,3}^{-1}+G_{\tau,3})\,dz =& 0
	\end{split}
\end{align}
Under the symmetric assumptions $p_2=3/4$ and $p_1+p_3=3/2$, the lhs of the
equations are actually the same.

\begin{figure}[hbt]
	\centering \includegraphics[width=\textwidth]{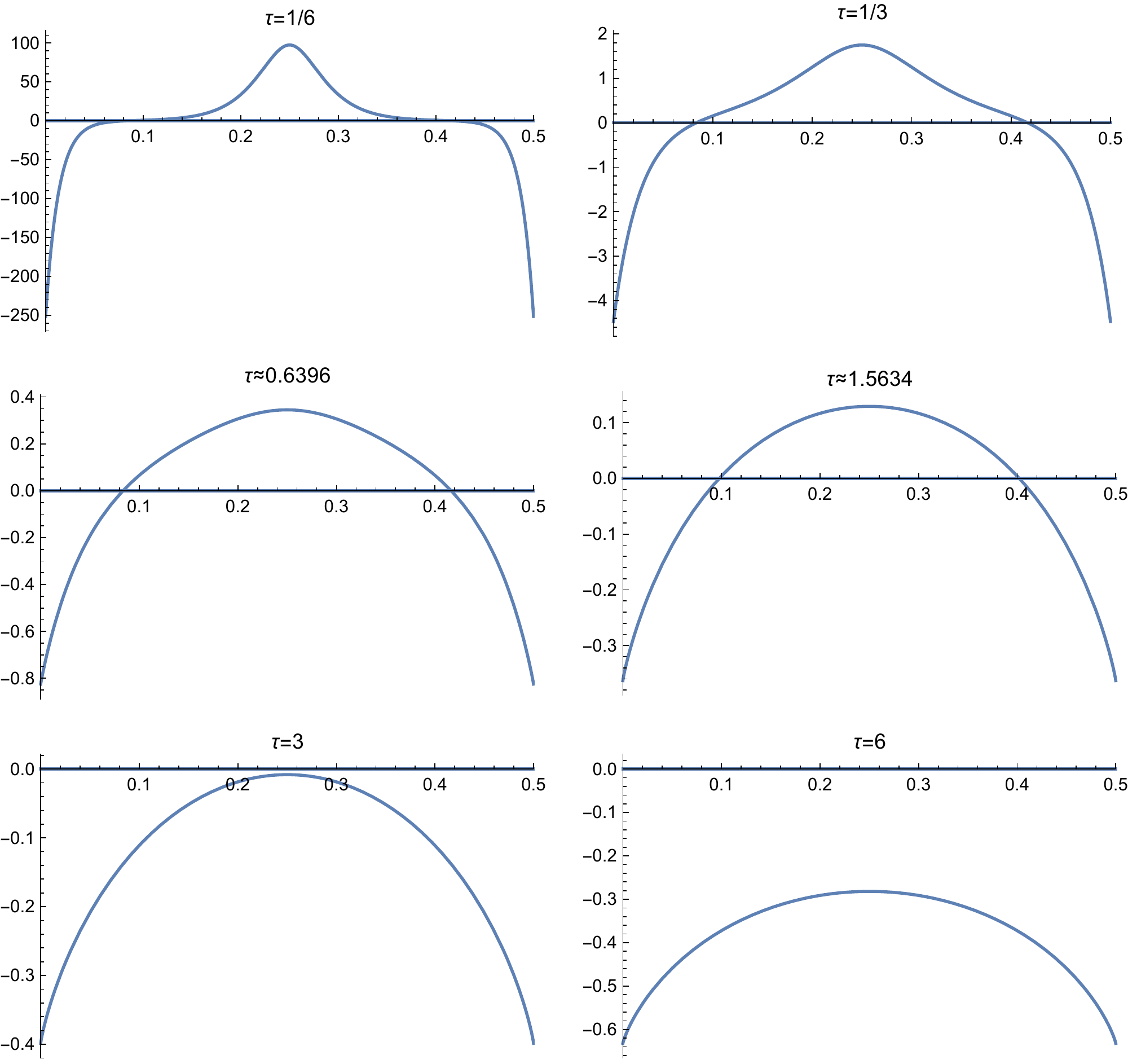}
	\caption{
		\label{fig:noroot} Plot of the lhs of~\eqref{eqn:period} against $p_1/3$
		for six different $\tau$'s.  In the middle row, the left is the twin of
		$\D$ surface, the right is the twin of $\P$ surface.
	}
\end{figure}

When we change $\tau$ to a very large value, Mathematica can not find any root.
In Figure~\ref{fig:noroot} we plot the lhs of~\eqref{eqn:period} against
$p_1/3$ for six values of $\tau$, which clearly shows absence of solution for
$\tau \ge 3$.  In these plots we allow $p_1>p_2=3/4$, in which case $p_1$
should be understood as $p_3$; the horizontal reflectional symmetry is then
obvious.  One observes that, as $\tau$ increases, $p_1$ and $p_3$ are pushed
towards the center $p_2=3/4$, finally meet there and vanish.   We calculate
that $p_1=p_2=p_3=3/4$ when $\tau=\tau_*=2.916517\cdots$.

\begin{observation}
	The period problem for $\Sigma_{\tau,3}$ has no solution for $\tau > \tau_*$.
\end{observation}

This observation is no surprise.  The smaller the parameter $t$ is, the closer
the catenoid unit is to a standard catenoid.   A standard catenoid has a
natural horizontal symmetry plane, hence it is easier to twin $S_t$ for small
$t$.  In the other limit $t \to \infty$, helicoids are forming.  As the only
properly embedded simply connected non-planar minimal surface in
$\RR^3$~\cite{meeks2005}, it is not possible to twin the standard helicoid.
Therefore, twinning $S_t$ with large $t$ is expected to be difficult.

Another way to understand the situation is the following: On the one hand, when
$t$ is large, $S_t$ is far from perpendicular to the horizontal planes.  To
meet the free boundary condition on the twin boundaries, a larger perturbation
would be needed.  Under the perturbation, the flat points at height $p_1$ and
$p_3$ is moving vertically towards $p_2$.  On the other hand, the flat points
are already very close in the $z$ direction, and there is little space for
perturbation.  These two factors together reduce the possibility of twinning
$\rPD$ surfaces with large $t$.

\subsection{Examples of small \texorpdfstring{$t$}{\it t}} \label{ssec:Traizet}

In the other end, for any $\delta$, the period problem is solved by
\cite{traizet2008} for sufficiently small $t$.  We now introduce his work.

Traizet was interested in TPMS with periods $(T_1,0)$, $(T_2,0)$ and
$(T_3,\varepsilon)$ which looks like infinitely many horizontal planes $H_0,
H_1, \cdots$ with $H_0$ at $z=0$ and $H_N$ at $z=\varepsilon$, and adjacent
planes are connected by small catenoid necks.  As $\varepsilon$ tends to $0$,
such a TPMS converges to a countably sheeted plane, and the catenoid necks
converges to singular points arranged in the lattice spanned by $T_1$ and
$T_2$.  We identify the plane to $\CC$ so that we can talk about positions and
periods by complex numbers.  Let $m_k$, $0 \le k < N$, be the number of
catenoid necks between $H_k$ and $H_{k+1}$ and $p_{k,i}$, $1 \le i \le m_k$, be
the limit position of the $i$-th catenoid neck.  The collection of
$\{p_{k,i}\}$ together with the periods $T_1$, $T_2$, $T_3$ is called a
\emph{configuration}.

Given a configuration, \cite{traizet2008} defines the \emph{force} on the
$i$-th neck between $H_k$ and $H_{k+1}$ as follows:
\begin{equation}\label{eq:balance}
	\begin{split}
		F_{k,i}=
		&\sum_{j \ne i}\frac{2}{m_k^2}\zeta(p_{k,i}-p_{k,j})\\
		&-\sum_{k'=k\pm 1}\sum_j \frac{1}{m_km_{k'}}\zeta(p_{k,i}-p_{k',j})\\
		&+\frac{1}{m_k}[(2x_k-x_{k-1}-x_{k+1})\eta_1+(2y_k-y_{k-1}-y_{k+1})\eta_2],
	\end{split}
\end{equation}
where
\[
	\zeta(z)=\frac{1}{z}+\sum_{0\ne w \in \left<
	T_1,T_2\right>}\Big(\frac{1}{z-w}+\frac{1}{w}+\frac{z}{w^2}\Big)
\]
is the Weierstrass $\zeta$ function; $\eta_i = 2\zeta(T_i/2)$, $i=1,2$; and the
center of mass $\sum_j p_{k,j}/m_k = x_kT_1 + y_kT_2$.  The configuration is
said to be \emph{balanced} if all the forces $F_{k,i}$ vanish, and
\emph{non-degenerate} if the differential of the map sending the positions
$\{p_{k,i}\}$ to the forces $\{F_{k,i}\}$ has real co-rank $2$.
\cite{traizet2008} proved that, if the configuration is balanced and
non-degenerate, then the TPMS described above exists and form a smooth family
for sufficiently small $\epsilon$.

The $\rPD$ surfaces with small $t$ are examples of Traizet's surfaces; the
configuration is given by
\begin{gather*}
	N=2,\quad m_0=m_1=1,\\
	T_1=1,\quad T_2=a=\exp(i\pi/3),\quad T_3=2(1+a)/3,\\
	p_{0,1}=0,\quad p_{1,1}=(1+a)/3.
\end{gather*}
The balance and non-degeneracy has been proved in Section 4.3.3 of
\cite{traizet2008}; see also Proposition~3 of the same paper.

\begin{figure}[hbt]
	\includegraphics[width=.3\textwidth]{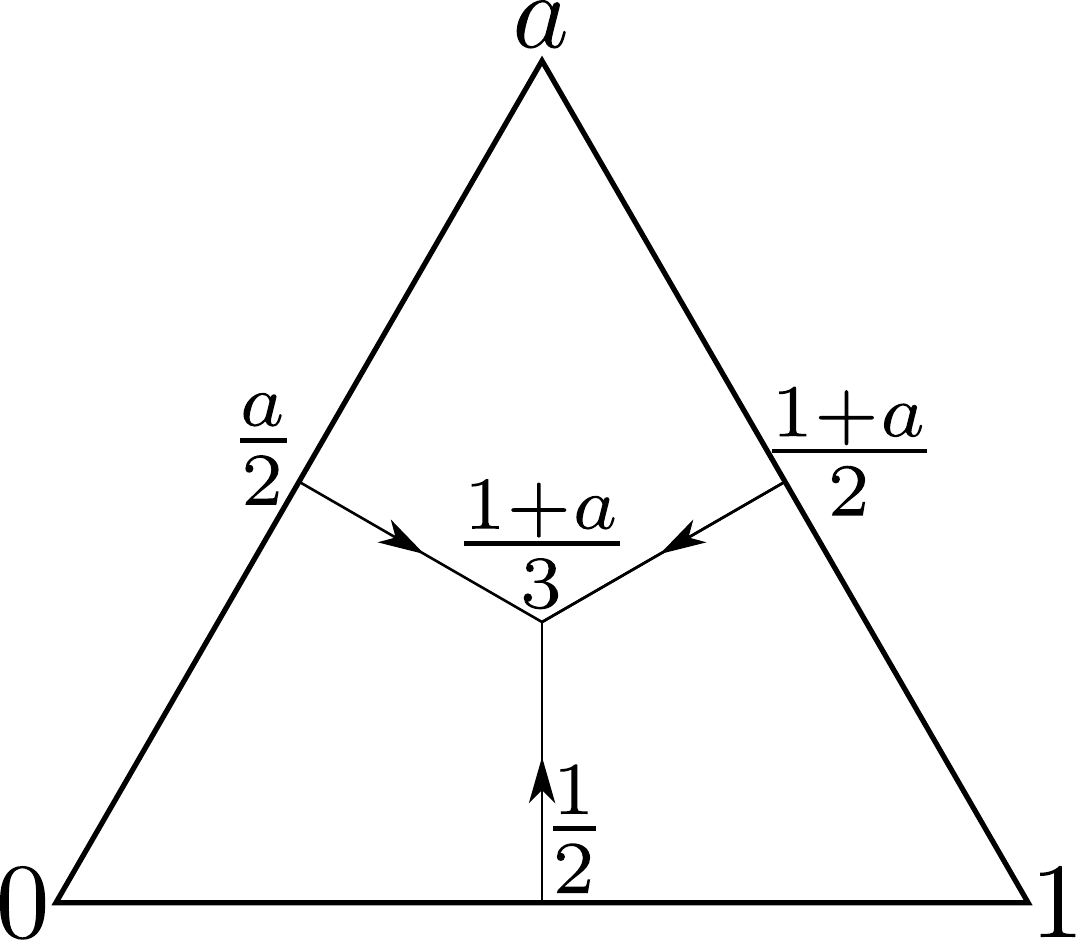}
	\caption{
		\label{fig:balance} The integrals of $\zeta'(z)=-\wp(z)$ along the three
		directed segments sum up to $0$ because of the symmetry.  This proves the
		force balancing condition for the $\H$ surfaces.
	}
\end{figure}

The $\H$ surfaces are also examples; the only difference from the $\rPD$
surfaces is $T_3=0$.  To verify Traizet's force balancing condition, consider
the three integral of $\zeta'(z)=-\wp(z)$ along the directed segments shown in
Figure~\ref{fig:balance}.  The Weierstrass elliptic function $\wp(z)$ is even
and invariant under the action of $\left<1,a\right>$.  The integrals sum up to
$0$ because of the symmetries.  That is,
\[
	\Big[\zeta\Big(\frac{1+a}{3}\Big) - \zeta\Big(\frac{1}{2}\Big)\Big] +
	\Big[\zeta\Big(\frac{1+a}{3}\Big) - \zeta\Big(\frac{a}{2}\Big)\Big] +
	\Big[\zeta\Big(\frac{1+a}{3}\Big) - \zeta\Big(\frac{1+a}{2}\Big)\Big] = 0.
\]
Since $\zeta(\frac{1}{2}) + \zeta(\frac{a}{2}) = \zeta(\frac{1+a}{2})$, we
conclude that
\[
	- 2\zeta\Big(\frac{1+a}{3}\Big) + \frac{4}{3}\zeta\Big(\frac{1}{2}\Big) +
	\frac{4}{3}\zeta\Big(\frac{a}{2}\Big)  = 0.
\]
Hence the configuration is balanced.  The non-degeneracy is verified
numerically.

At this point, crystallographers should notice the analogy of $\rPD$ and
$\H$ family with the cubic and hexagonal closed packings.  Within a period,
there are three possible limit positions for the necks, namely $j(1+a)/3$
with $j=0, 1, 2$.  We then obtain an infinite sequence consisting of letters
$0$, $1$ and $2$.  The $\rPD$ surfaces has a periodic sequence
$\cdots012012012\cdots$.  The sequences for $\rPD$ twins should be periodic and
palindromic, like $\cdots010101010\cdots$ for the $\H$ surface,
$\cdots012101210\cdots$ for the T-WP surface, etc.

The formula~\eqref{eq:balance} only involves adjacent layers.  Up to symmetry,
the only possible consecutive triples is $012$ and $010$, for which the force
balancing condition has been verified in $\rPD$ and $\H$ surfaces.  So the
$\rPD$ twin configurations are all balanced and non-degenerate.  Consequently,
for any $\delta>0$, there is a real number $\epsilon_\delta>0$ such that the
minimal twin surface $\Sigma_{\tau,\delta}$ exists for
$0<\tau<\epsilon_\delta$.  In fact, the argument works for any periodic
sequence. For example, the sequence $\cdots0101201012\cdots$ also satisfy the
force balancing condition.

Traizet's proof uses implicit function theorem.  We do not know if the
following conjecture holds.
\begin{conjecture}\label{conj:approx}
	There is a positive number $\epsilon>0$ such that the $\rPD$ twin surface
	$\Sigma_{\tau,\delta}$ exists for all $0 < \tau < \epsilon$ and $\delta>0$.
\end{conjecture}
A proof of this conjecture would justify our strategy of using polysynthetic
twins to approximate the single reflection twin.

\subsection{Numerical examples with Surface Evolver} \label{ssec:rPDevolver}

Traizet's result guarantees TPMS near the catenoid limit.  It does not tell how
small is sufficient.  Sometimes $t$ doesn't need to be very small.  In the case
of $\rPD$ and $\H$ surfaces, the Traizet families extend all the way to the
helicoid limit.  For a concrete surface alleged in a Traizet family, we need
numerical technique to confirm its existence.

Brakke's Surface Evolver~\cite{brakke1992} is a software that simulates the
physics of surfaces by minimizing energies.  Surfaces are modeled by
triangulations, and the energy is minimized by gradient descent method.
Surface Evolver can handle various energies under various constraints.

In a usual application, the energy to minimize is the surface tension energy or
the area functional.  However, contrary to the intuition of many, TPMS do not
result from area minimization.  In the translational unit, a TPMS is actually a
strict maximum of the area functional among its parallel
surfaces~\cite{kgb2012}.  A TPMS is indeed stable for the area functional if a
volume constraint is imposed on the translational unit~\cite{kgb1996}.  But
this is again false for any slightly larger piece~\cite{ross1992}.  Between two
adjacent twin boundaries, $\Sigma_{t,\delta}$ correspond to a large piece of
$S_t$, hence the area functional is not our option.

We will minimize instead the Willmore energy, or more precisely, the integral
of the squared mean curvature; see~\cite{hsu1992} for Surface Evolver
experiments on this energy.  Physically, the Willmore energy measures the
deviation of the surface from the zero mean curvature, thus models the
``bending'' energy.  The Willmore energy vanishes on minimal surfaces.

We construct the $\rPD$ twin surface $\Sigma_{t,\delta}$ in Surface Evolver in
four steps:

\begin{enumerate}
	\item Prepare a catenoid unit.

		Let $h=h(t)$ be the height function defined in Section~\ref{ssec:rPD}.
		Consider four points
		\[
			A(1,0,0), B(-1/2,\sqrt{3}/2, 0), C(-1,0,h), D(1/2,\sqrt{3}/2,h).
		\]
		The catenoid unit of $S_t$ modulo reflectional symmetries is a minimal
		surface with fixed boundary conditions on the segments $AB$ and $CD$, and
		free boundary conditions on the vertical planes $y \pm \sqrt{3} x =
		\sqrt{3}$.  Reflections in these planes yield the whole catenoid unit.  The
		boundary triangles has unit inradius, compatible with the discussions in
		Section~\ref{ssec:rPD}.  The unit is obtained in Surface Evolver by
		minimizing the area functional.  As a consequence of this practice, we can
		only obtain $\rPD$ surfaces with parameter $t>t_0=0.494722\cdots$.

	\item Generate a slab.

		We generate a slab of $S_t$ from the catenoid unit by order-$2$
		rotations about the segment $AB$ or $CD$.  In Surface Evolver, this is done
		by listing the matrices, say \texttt{a} and \texttt{b}, of the two rotations in
		\texttt{VIEW\_TRANSFORM\_GENERATORS}.  Assume an odd $\delta$, we set
		\texttt{transform\_expr} to be a string $k$\texttt{(ba)} if $\delta=4k-1$ or
		$k$\texttt{(ba)b} if $\delta=4k+1$.  Then the Surface Evolver will display $\delta+1$
		catenoid unit.  Finally, we use the \texttt{detorus} command to convert the
		displayed slab unit into a real surface.

	\item Slice the slab.

		In this step, we use Brakke's script \texttt{slicer.cmd}.  It is included
		in the Evolver distribution, and also available on the website of Surface
		Evolver.  It removes the part of the surface on one side of a given plane.
		We slice the slab in the previous step by two horizontal planes with offset
		$0.5$ at distance $\delta$ apart.  The script also marks the vertices and
		edges newly created by the slice.  This allows us to impose free boundary
		conditions by constraining the new vertices and edges on the slicing
		planes.

	\item Evolve the surface.

		We turn off the surface tension energy (\texttt{set facet tension 0}), and
		turn on the Willmore energy, which is the integral of squared mean
		curvature.  Then we leave Evolver to minimize the energy.  Apart from the
		command \texttt{g} that does one step of gradient descent, the command
		\texttt{hessian\_seek} is particularly useful in this step to accelerate the
		calculation.  If the Willmore energy decreases to practically $0$, we
		obtain a minimal surface with free boundary condition on the boundary of a
		triangular prism.  Reflections in the faces of the prism give the whole
		$\Sigma_{t,\delta}$.
\end{enumerate}

\begin{figure}[hbt]
	\includegraphics[width=.7\textwidth]{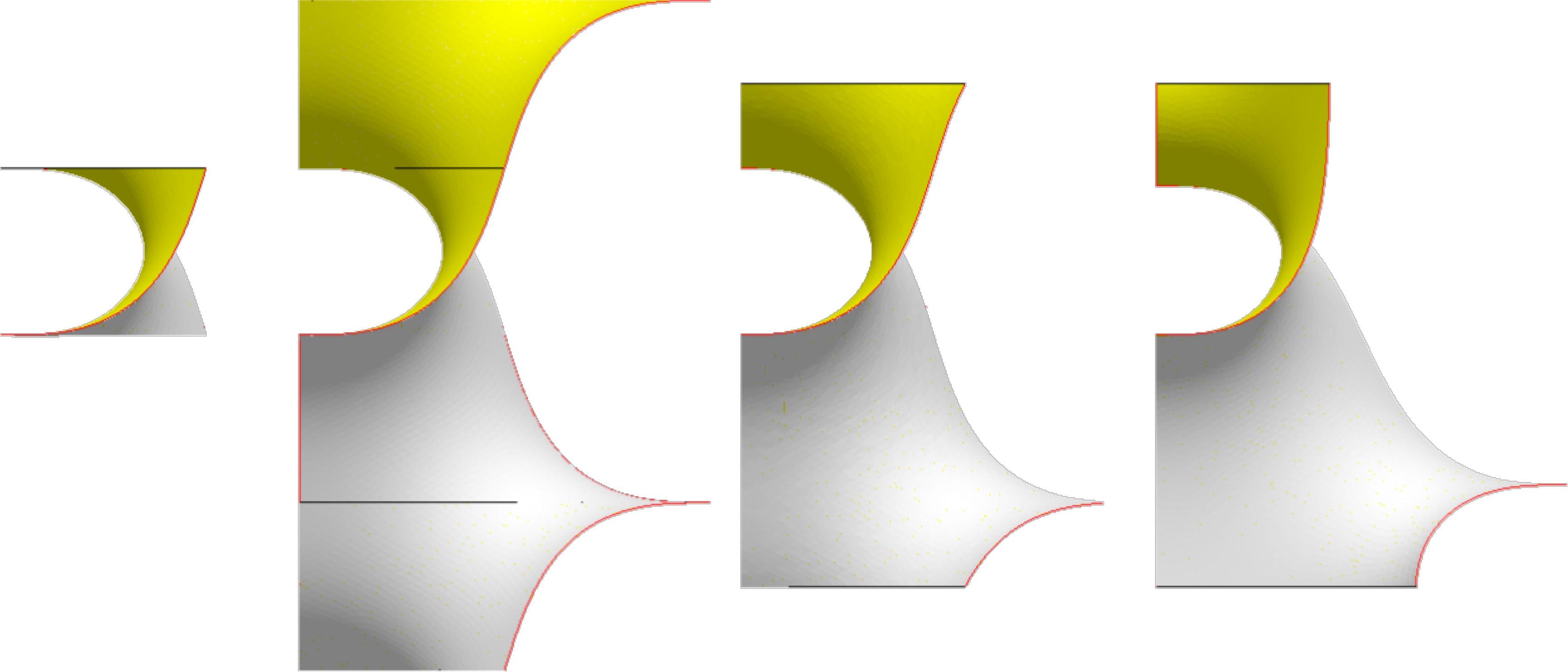}
	\caption{\label{fig:evolver} Result of each step with $t=\sqrt{2}$ and $\delta=3$.}
\end{figure}

In Figure~\ref{fig:evolver}, we show the result of each step with $t=\sqrt{2}$
(Schwarz' $\P$ surface) and $\delta=3$.

The program can be easily modified to change the parameters $t$ and $\delta$.  We
perform calculations for parameters $t=t_0, \sqrt{1/2}$ ($\D$), $1, \sqrt{2}$
($\P$), and $\delta=1$ ($\H$), $3, 5, 11, 21, 99$.  When $\delta$ is big, the last step
will be very slow.  We can increase the speed in the price of precision by
reducing the number of faces in the catenoid unit.  In Surface Evolver, the
command \texttt{r} refines the surface by subdividing each triangle into four.
When generating the catenoid unit, we refine five times for $\delta=1,3,5$, four
times for $\delta=11,21$ and three times for $\delta \ge 99$.

\begin{figure}[bt]
	\includegraphics[width=.7\textwidth]{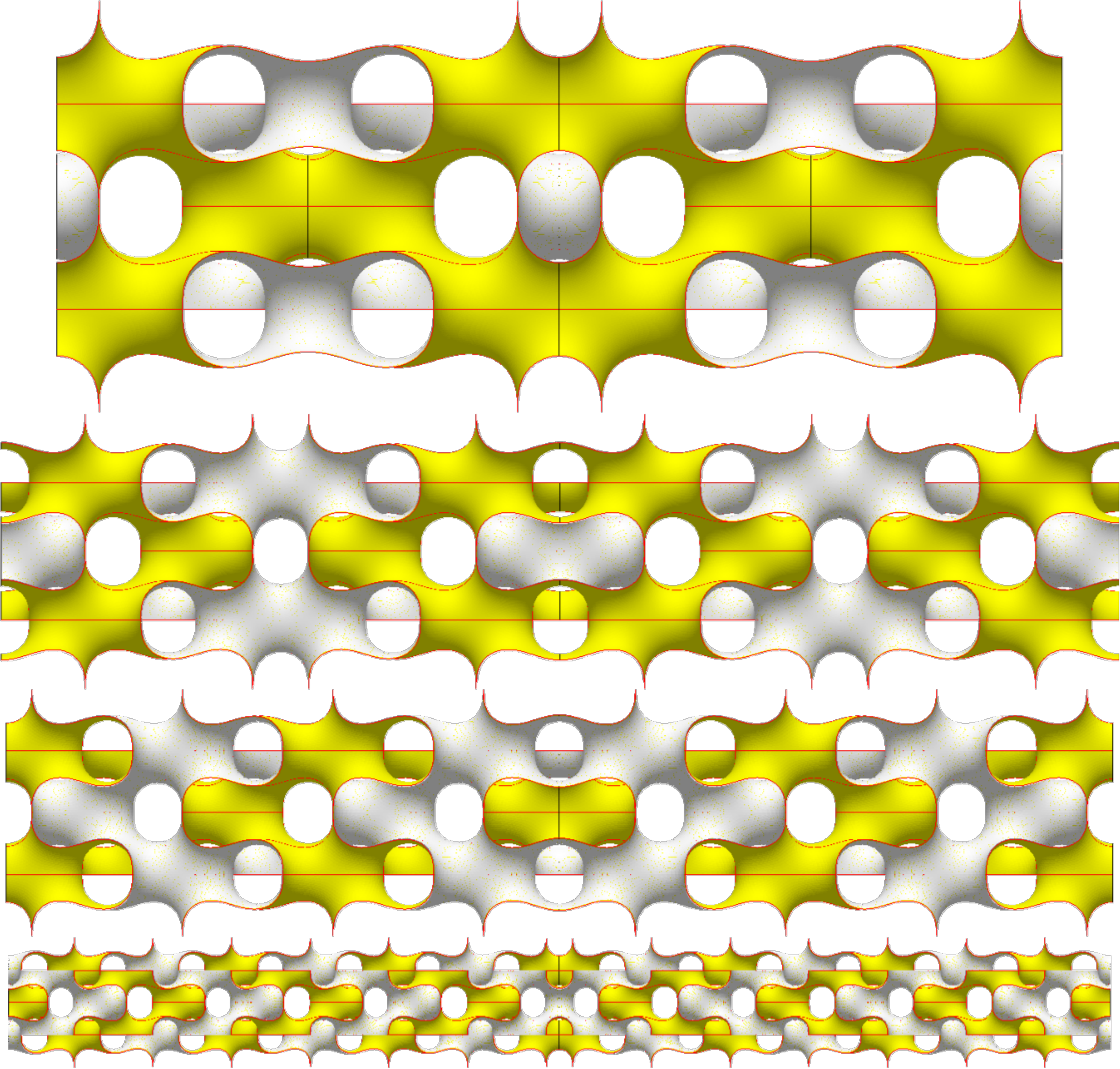}
	\caption{
		\label{fig:resultD} $\D$ twin surfaces ($t=\sqrt{1/2}$).  From top to bottom,
		$\delta=3, 5, 11, 21$.  The twin boundary is clearly seen in the middle.
	}
\end{figure}

\begin{figure}[bt]
	\includegraphics[width=.7\textwidth]{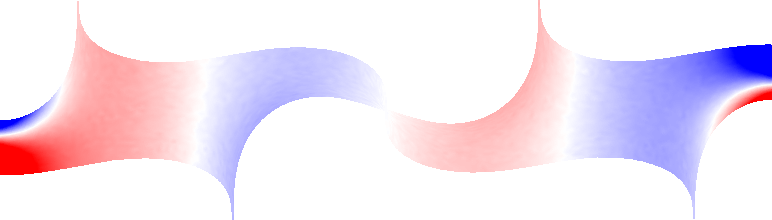}
	\caption{
		\label{fig:deviationD} Deviation of the $\D$ twin with $\delta=5$ from the
		$\D$ surface.  The picture shows the slab between two twin boundaries.
		Blue and red indicate perturbations in opposite directions; white indicates
		coincidence.
	}
\end{figure}

The calculation goes surprisingly smooth for small $t$ and small $\delta$.  For
$t=t_0$ and $t=\sqrt{1/2}$ ($\D$), the Willmore energy decreases quickly to the
order of $10^{-26}$ or lower, even when $\delta=99$.  In view of the discussion
in Section~\ref{ssec:Weber}, it is reasonable to believe that the calculation
would be even faster for smaller $t$.  The results for $t=\sqrt{1/2}$ and
$\delta=3,5,11,21$ are shown in Figure~\ref{fig:resultD}.  The difference from
the $\D$ surface is not visible with human eyes.  We use
CloudCompare~\cite{cloudcompare} to calculate the deviation from the $\D$
surface; the result for $\delta=5$ is colored accordingly in
Figure~\ref{fig:deviationD}.  We can see that the perturbation decays away from
the twin boundaries, and the flat points are perturbed towards the middle.
These pictures justify that the structures observed in~\cite{han2011} are
indeed minimal $\D$ twins.

\begin{observation}
	The $\rPD$ twin surfaces $\Sigma_{t,\delta}$ exist for a large set of parameters.
\end{observation}

In fact, if we obtain a minimal surface, the calculation roughly solves the
period problem.  More specifically, the flat points are the vertices
transformed from $A$ or $C$ when we generate the slab.  They correspond to the
poles and zeros of the Gauss map $G_{\tau,\delta}$.  After normalization, the
heights of these points give the sequence $(p_k)_{1 \le k \le \delta}$.  For
example, the $\Sigma_{\sqrt{1/2}, 3}$ ($\D$ twin) generated by Surface Evolver
has twin boundaries at $z=\pm 3\sqrt{2}/2$, and the surface intersects the
lateral edges of the prism at $z \approx 0$ and $z \approx \pm 1.41240$.  By a
scaling and a translation, we may place the twin boundaries at $z=0$ and
$z=3/2$, and obtain $p_1 \approx 0.25064$, $p_2 \approx 0.75$ and $p_3 \approx
1.24936$.  Compared to the solution given by the Mathematica program used
for~\cite{fujimori2009}, these values are accurate until the fifth digit after
the decimal point.  However, the precision is not guaranteed as $t$ increases.
The $\Sigma_{\sqrt{2}, 3}$ ($\P$ twin) shown in Figure~\ref{fig:evolver} gives
$p_1 \approx 0.303350$, $p_2 \approx 0.75$ and $p_3 \approx 1.196651$, while
Mathematica computes more precisely $p_1 = 0.293406\cdots$, $p_2=0.75$ and $p_3
= 1.20659\cdots$.

\medskip

The Willmore energy also decreases to the order of $10^{-26}$ for $t=1$ with $\delta
\le 21$, and for $t=\sqrt{2}$ with $\delta\le 3$.  However, for $t=1$ and $\delta=99$,
the Surface Evolver only manages to reduce the Willmore energy to the order of
$10^{-8}$; then the decrement becomes extremely slow.  For $t=\sqrt{2}$
($\P$ surface) and $\delta=5$, the energy seems to stabilize at the order $10^{-4}$.
These calculations are therefore inconclusive.

\begin{observation}
	Surface Evolver does not converge to a minimal surface in a reasonable time
	if $t$ and $\delta$ is too large.
\end{observation}

We then append more calculations for $t=\sqrt{1/2}$ to see if it eventually
fails with sufficiently large $\delta$.  It turns out that we are able to
obtain a minimal surface with $\delta \le 299$.  Since the calculation is very
time-consuming for large $\delta$, we did not perform further computation.  But
the success with large $\D$ twins ensures that the inconclusive results for
$\P$ twins do not arise from numerical or system errors.

\begin{remark}
	Despite the negative evidences, we refrain from conjecturing that $\P$ twins
	do not exist.  Recently, \cite{han2017} observed interconversions between
	$\P$ surface and $\D$ surface in their experiment.  The structure looks like
	a $\P$ surface and a $\D$ surface glued along a $(111)$ lattice plane.  It
	suggests a possible mechanism of twinning the $\P$ surface through a $\D$
	twin.  However, it is too complicated to calculate the Weierstrass
	representation for large $\delta$, and Surface Evolver might be trapped in a
	local minimum of the Willmore energy, therefore won't reach such a
	configuration.  So the approaches of this note wouldn't reveal such a
	twinning.

	\cite{morabito2012} and~\cite{traizet2013} connected catenoid necks between
	countably many horizontal planes in a non-periodic way.  In their works,
	every adjacent pair are connected by finitely many necks.  But it is possible
	to extend their technique to connect necks between countably many tori in a
	non-periodic way\footnotemark; see~\cite{traizet2008}.  This would lead to a
	rigorous construction for reflection $\rPD$ twins, as well as the $\P$--$\D$
	interconversion and many other interesting structures, and would not involve
	any approximation as we hoped in Conjecture~\ref{conj:approx}.  This will be
	the topic of a future project.
\end{remark}

\footnotetext{Confirmed by Traizet through personal communication.}

\section{G twin surfaces and some speculations} \label{sec:GTwin}

\subsection{Numerical examples with Surface Evolver} \label{ssec:Gevolver}

Twinnings of the $\G$ surface are also observed in experiment.  Taking the flat
points as lattice points, the observed twin boundaries are $\{211\}$ planes
with offset $0$.  

To generate $\G$ twin surfaces, we follow the same procedure for $\rPD$
twinning:  Take a TPMS and impose free boundary condition on two parallel
planes.  If Surface Evolver is able to reduce the Willmore energy to
practically $0$, then we obtain a minimal twin surface.  A few modifications is
however necessary.

\begin{enumerate}
	\item Prepare the $\G$ surface in an orthorhombic cell.

		Unlike the $\rPD$ surfaces, the $\G$ surface contains no straight line, and
		has no symmetry plane.  A datafile of initial triangulation is prepared by
		Fujita using the \texttt{torus} model of Surface Evolver.  The orthorhombic
		unit cell is generated by three vectors in the $[0\bar{1}1]$, $[\bar{1}11]$
		and $[211]$ directions.  For our convenience, we made a small modification
		to place the $[211]$ direction on the $z$-axis.

		Since this unit cell is larger than the translational fundamental cell, it
		is not safe to minimize the area functional~\cite{kgb1996}.  We prepare the
		$\G$ surface by minimizing the Willmore energy directly.

	\item Generate a slab.

		In the torus model, the command \texttt{y 3} duplicates the displayed
		surface along the $z$-axis ($[211]$-direction).  Repetition of the command
		for $k$ times would generate $2^k$ copies of the orthorhombic cell.  We do
		not \texttt{detorus} the surface since we need the double periodicity.

	\item Slice the slab.

		Brakke's script does not work well in the torus model.  We write our own
		script to slice the slab by two horizontal planes ($(211)$ planes) at
		offset $0$, and impose free boundary conditions on the slicing planes.

	\item Evolve the surface.
\end{enumerate}

\begin{figure}[hbt]
	\includegraphics[width=0.7\textwidth]{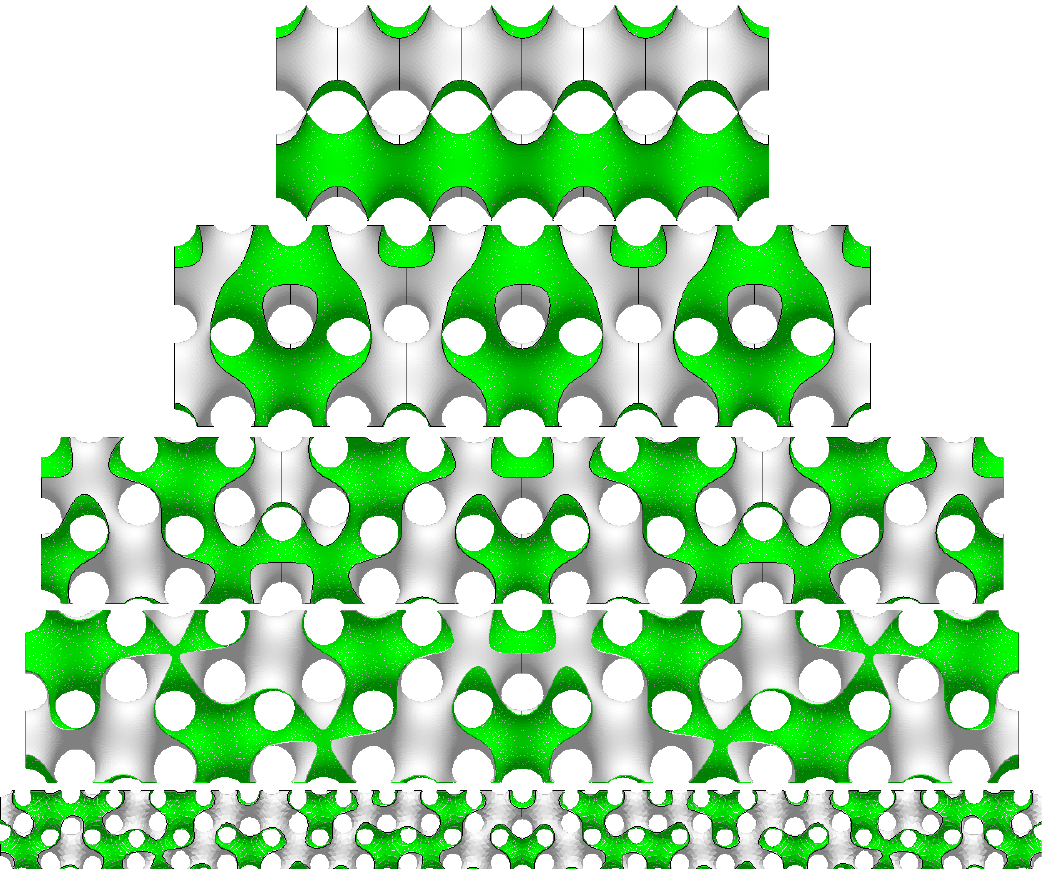}
	\caption{
		\label{fig:resultG} $\G$ twin surfaces.  From top to bottom, $\delta=1, 2, 5,
		10, 23$.  The twin boundary is clearly seen in the middle.
	}
\end{figure}

\begin{figure}[hbt]
	\includegraphics[width=0.5\textwidth]{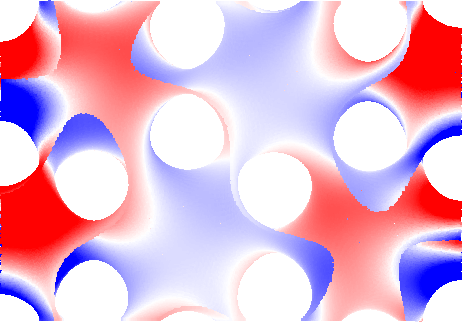}
	\caption{
		\label{fig:deviationG} Deviation of the $\G$ twin with $\delta=5$ from the
		$\G$ surface.  The picture shows the slab between two twin boundaries.
		Blue and red indicate perturbations in opposite directions; white indicates
		coincidence.
	}
\end{figure}

The program can be easily modified to adjust the lattice distance $\delta$
between the slicing planes (twin boundaries).  We perform calculation for
$\delta=1, 2, 5, 10, 23$.  In all these cases, we are able to obtain a minimal
surface by reducing the Willmore energy to the order of $10^{-26}$.  We did not
try larger $\delta$.  The results are shown in Figure~\ref{fig:resultG}.  The
$\G$ twin with $\delta=1$ turns out to be an orthorhombic deformation of the
$\D$ surface (an oD surface); this will be explained later.  In
Figure~\ref{fig:deviationG}, the result for $\delta=5$ is colored according to
the deviation from the $\G$ surface.

\begin{observation}
	The $\G$ twin surfaces exist for a large set of parameters.
\end{observation}

\subsection{Cubic polyhedral models} \label{ssec:cubic}

The Gyroid surface is probably the most beautiful TPMS, but has the reputation
of being difficult to visualize.  For Surface Evolver, a datafile written by
Gro{\ss}e-Brauckmann back in 1995 is still widely used to generate the $\G$
surface.  It decomposes the cubic unit cell of the $\G$ surface into hexagons
and rectangles, which is difficult to write from scratch.

\begin{figure}[hbt]
	\includegraphics[width=0.7\textwidth]{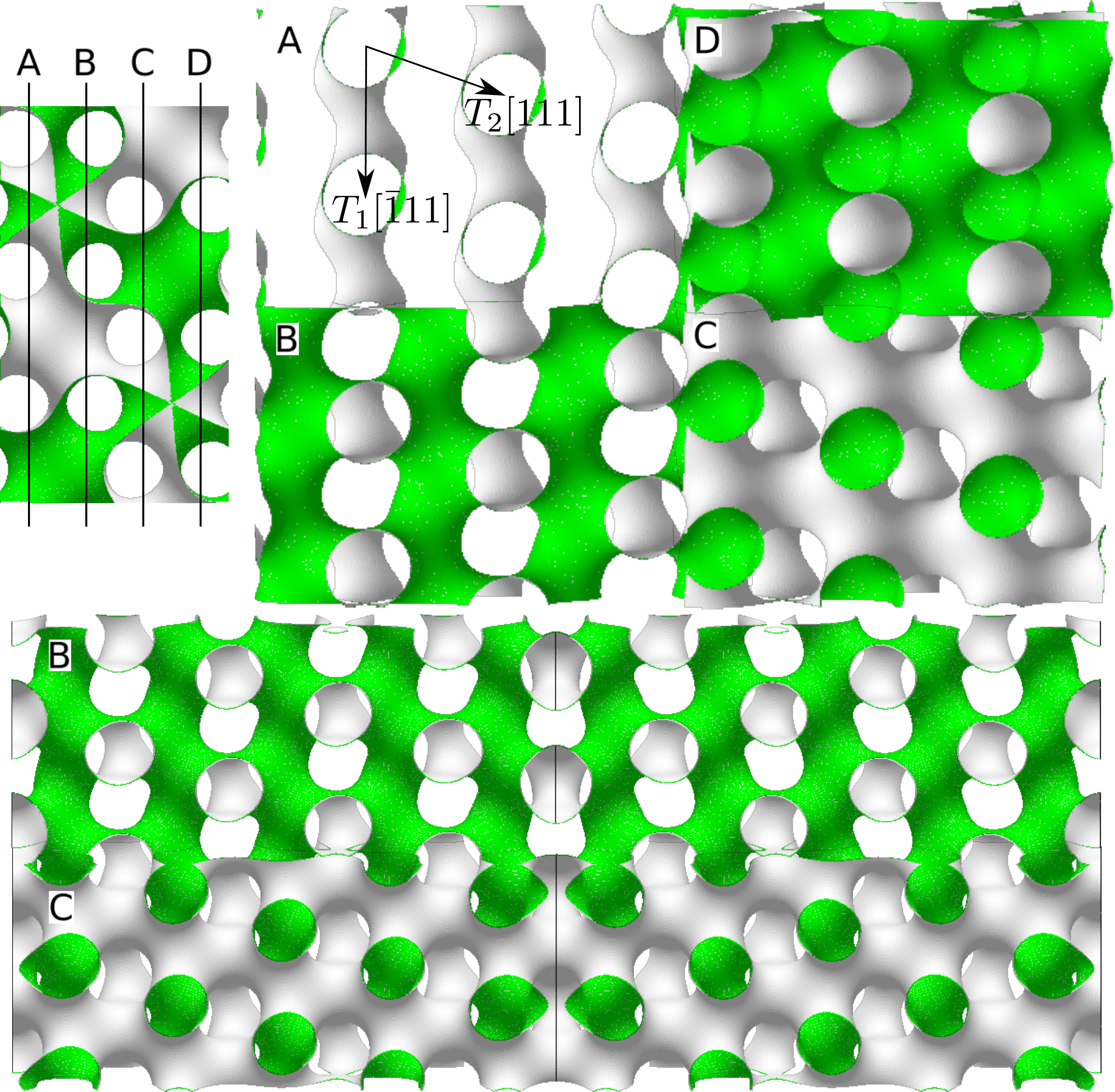}
	\caption{
		\label{fig:G110} The $\G$ surface (upper right) and the $\G$ twin with
		$\delta=10$ (bottom) sliced by $(0\bar{1}1)$ planes at offset $0.5$.  The
		position of the slicing planes are indicated in the upper left corner,
		where the $\G$ surface is seen in the $[\bar{1}11]$ direction.  In the
		upper right corner, intersection of the $\G$ surface with the slicing
		planes is highlighted by blue cycles.
	}
\end{figure}

The $(211)$ $\G$ twins reveal an interesting structure of the $\G$ surface
which, to the knowledge of the author, was not mentioned before.  In the upper
half of Figure~\ref{fig:G110} is a gyroid sliced by $(0\bar{1}1)$ planes at
offset $0.5$.  We observe necks arranged in a $2$-D lattice spanned by vectors
$T_1[\bar111]$ and $T_2[111]$, and the lattices in adjacent layers differ by
$T_1/2$ or $T_2/2$, alternatively.  In the lower half of the same figure, we
show the slices of the $\G$ twin with $\delta=10$.  We see clearly that the
twin boundary is parallel to $T_1$ and perpendicular to the $(0\bar{1}1)$
planes.

\begin{figure}[hbt]
	\includegraphics[width=0.6\textwidth]{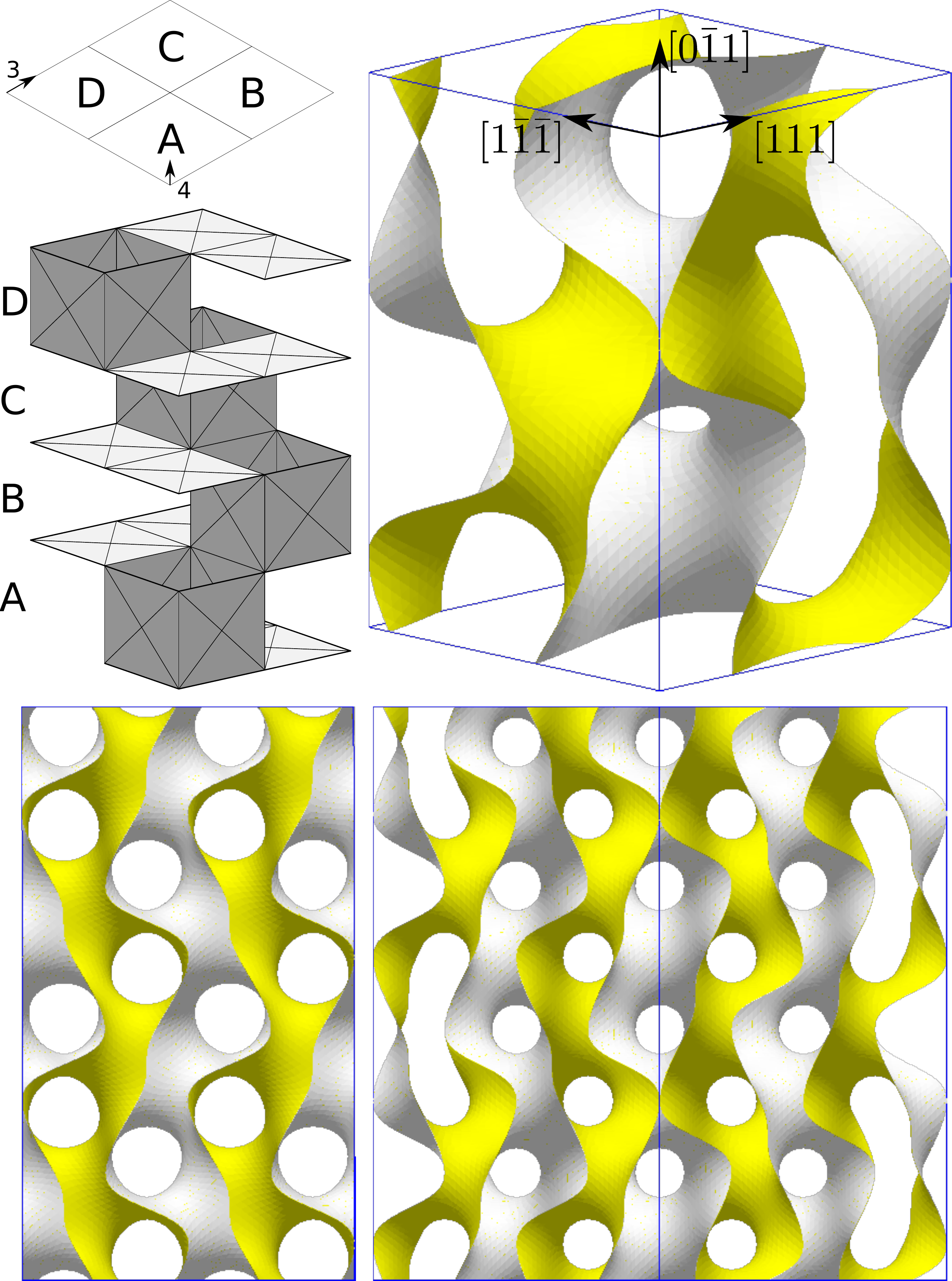}
	\caption{
		\label{fig:NewG} A new construction of the Gyroid in Surface Evolver.
		Upper left: the initial surface consists of four horizontal planes
		connected by tubes; see text.  The rhombic tubes are labeled in accordance
		with slicing planes in Figure~\ref{fig:G110}.  The small arrows indicate
		the directions of $3$- and $4$-fold axis.  Upper right: Evolved result in a
		monoclinic unit cell.  Lower left: Evolved result seen along the $3$-fold
		axis.  Lower right: Evolved result seen along the $4$-fold axis.
	}
\end{figure}

This observation leads to a new way to construct the $\G$ surface in Surface
Evolver.  The initial triangulation in torus model is shown in the upper left
corner of Figure~\ref{fig:NewG}.  The monoclinic unit cell is spanned by
$(\pm\sqrt{2}, 1, 0)$ and $(0,0,2\sqrt{2})$; this cell has the same volume as
the cubic cell.  We prepare four horizontal planes at $z=j\sqrt{1/2}$,
$j=0,1,2,3$.  Each plane is a $2$-torus decomposed into four rhombi.  Vertical
rhombic tubes connect adjacent planes as shown in the figure.  After several
refinements and evolvement, this simple configuration converges correctly to the
$\G$ surface by minimizing the Willmore energy; see the evolved result in
Figure~\ref{fig:NewG}.

Note that the initial triangulation is a discrete surfaces with quadrangle
faces, containing all the vertices and edges from the monoclinic lattice.  We
call such surfaces \emph{cubic polyhedral}.  For cubic lattice, cubic
polyhedral surfaces have been studied by~\cite{goodmanstrauss2003} as a
discrete model for Schwarz' $\P$ and CLP surfaces.  Hence we just extended
cubic polyhedral model to a monoclinic lattice for the $\G$ surface.

\medskip

\begin{figure}[hbt]
	\includegraphics[width=0.6\textwidth]{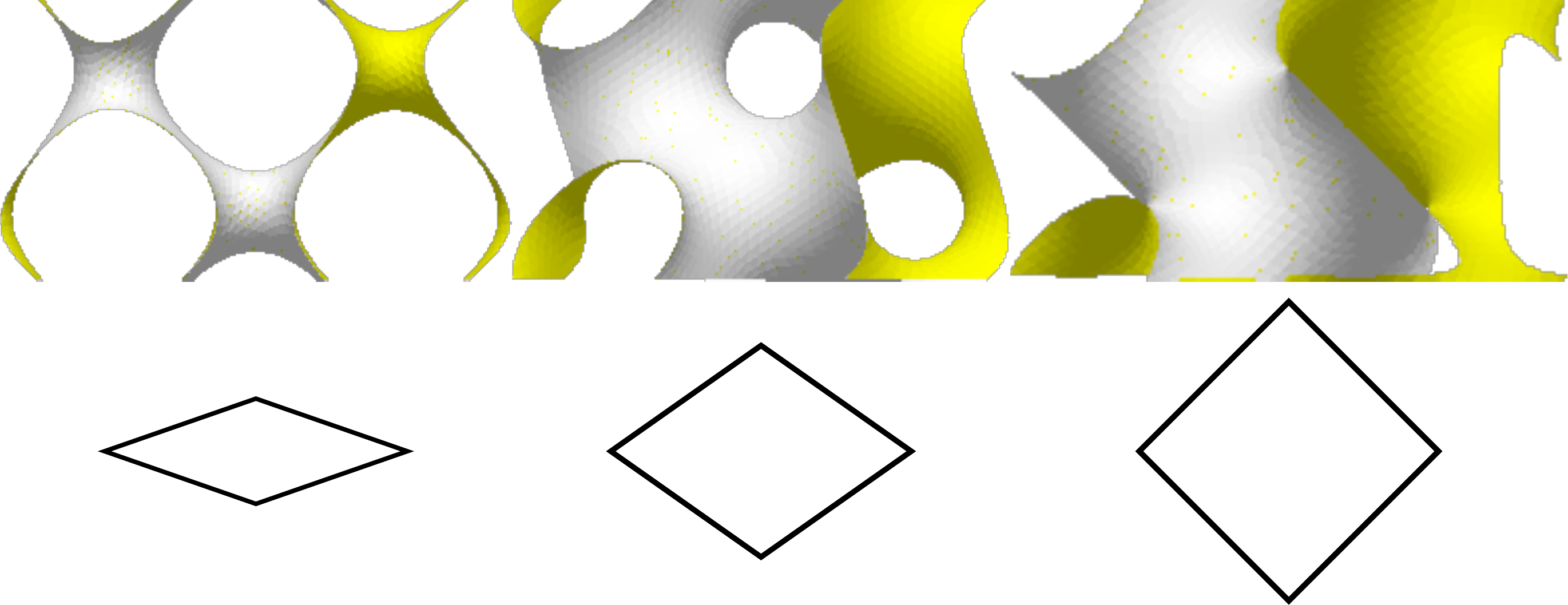}
	\caption{
		\label{fig:mG}  The diagonal of the rhombus (vertical in the figure) is a
		$4$-fold axis of the $\G$ surface.  By stretching the $\G$ surface (middle)
		in this direction, we obtain the $\D$ surface (right) and a $4$-fold saddle
		tower (left).  The surfaces shown in this figure are in the triclinic unit
		cell; the $4$-fold axis is towards the readers.  To be compared with $\tG$
		surfaces in Figure~11 of \cite{fogden1999}.
	}
\end{figure}

The rhombus spanned by $(\pm\sqrt{2}, 1, 0)$ can be stretched along the
diagonal (direction $(0,1,0)$, $[100]$ in the gyroid lattice) into a rhombus
spanned by $(\pm\sqrt{2}, c, 0)$.  It becomes a square of side length $2$ when
$c = \sqrt{2}$.  We then repeat the experiment with different values of $c$.
The surface deforms continuously, and becomes the $\D$ surface when $c =
\sqrt{2}$; see the right side of Figure~\ref{fig:mG}.  The horizontal planes
are $(100)$ planes of the diamond lattice.  Indeed, the diamond lattice has the
symmetry of square lattice on $(100)$ planes, and adjacent $(100)$ planes
differ by a half-period in alternating directions.  We then obtain a new cubic
polyhedral model for the $\D$ surface, and for every surface along this
deformation.

\begin{remark}
	In the experiments that produce Figure~\ref{fig:mG}, the initial
	triangulation consists of two horizontal planes a triclinic unit cell spanned
	by $(\pm\sqrt{2}, c, 0)$ and $(0, c, \sqrt{2})$.  Surface Evolver runs
	faster and better on this smaller unit cell.
\end{remark}

The deformation above provides a continuous path from $\G$ to $\D$.  This
seems strange because the $\D$ surface is a non-degenerate surface in the
five-dimensional Meeks family, hence should stay in the family after small
perturbations; but the $\G$ surface does not belong to Meeks
family\footnotemark.  Indeed, a closer look reveals that the deformation above
is not direct.  As we increase $c$, the deformation seems to follow first the
$\tG$ family (not in Meeks family), then the $\tD$ family (in Meeks family);
see~\cite{fogden1999} for details about $\tD$ and $\tG$ families.  We now
present some evidences.

Note that the $[100]$ direction is a $4$-fold axis of the $\G$ surface.  If we
reduce $c$, the 2D lattice is compressed into a very elongated rhombic lattice,
and we obtain a $4$-fold ``saddle tower''; see the left side of
Figure~\ref{fig:mG}.  This behavior is compatible with the $\tG$ family.  The
tunnels that we see through in Figure~\ref{fig:mG} correspond to helix curves
in the surfaces.  As we increase $c$, these tunnels shrink, since helix curves
are becoming straighter.  They eventually degenerate into straight lines as the
deformation enters the $\tD$ family.  Our experiment shows that the degeneracy
occurs before we reach the $\D$ surface.  We are not able to determine
precisely where the degeneracy takes place, but it seems to be around the tGD
surface as described in~\cite{fogden1999}.

\footnotetext{The author thanks Mathias Weber and the anonymous referee for
pointing this out.}

\begin{remark}
	\cite{sadoc1989} observed that, by contracting edges, one can transform the
	skeleton of the $\G$ surface, into the skeleton of the $\D$ surface, and
	finally into the skeleton of $\P$.  Between $\D$ and $\P$, such
	transformation is explicitly given by the $\rPD$ family.  Between $\G$ and
	$\D$, our experiment suggest that such transformation is provided by the
	$\tG$-$\tD$ path.
\end{remark}

\subsection{A monoclinic family of TPMS} \label{ssec:mG}

Inspired by the experiments above, we notice the following family of Traizet
configurations, which turns out to be balanced for any linearly independent
$T_1$ and $T_2$ in $\CC$, and is generically non-degenerate: 
\begin{gather} \label{eq:mG}
	\begin{split}
		&N=2,\quad m_0=m_1=1,\\
		&p_{0,1}=0,\quad p_{1,1}=T_1/2,\\
		&T_3=(T_1+T_2)/2.
	\end{split}
\end{gather}

It follows from \cite{traizet2008} that each configuration~\eqref{eq:mG}
implies a family of TPMS near the catenoid limit, whose horizontal lattice
planes admit the symmetry of the 2D lattice spanned by $T_1$ and $T_2$.  Near
Traizet's limit, we can also horizontally deform the rhombic lattice into other
2D lattices, as we did in the experiments above.  In particular, there is a
Traizet family ($T_{1,2}=\pm\sqrt{2}+i$) whose horizontal planes has the
symmetry of the $(011)$ planes of the $\G$ surface, and another Traizet family
($T_2 = i T_1$) whose horizontal planes has the symmetry of the $(100)$ planes
of the $\D$ surface.  We name these alleged deformations the $\mG$ family,
since the surfaces have the symmetry of monoclinic lattices.  In the
neighborhood of a generic surface, this would be a three dimensional family:
one dimension comes from the vertical deformations; the other two from the
horizontal deformations.

One would naturally expect that every surface along the deformation path
between the $\D$ and $\G$ surfaces belongs to a Traizet family.  This seems not
the case, at least not directly, for the $\D$ surface, as we now explain.  We
speculate the same phenomena for the $\G$ surface.

Note that the vertical direction ($[100]$ in the diamond lattice) is the
$4$-fold axis of the $\D$ surface.  Hence the Traizet family should preserve
the tetragonal symmetry.  We already know a tetragonal deformation family for
the $\D$ surface, namely the $\tD$ family~\cite{fogden1999}.  However, $\tD$
family admits a helicoid limit, but no catenoid limit.  This suggests that the
Traizet family and the $\tD$ family are not the same.

\begin{figure}[htb]
	\includegraphics[width=\textwidth]{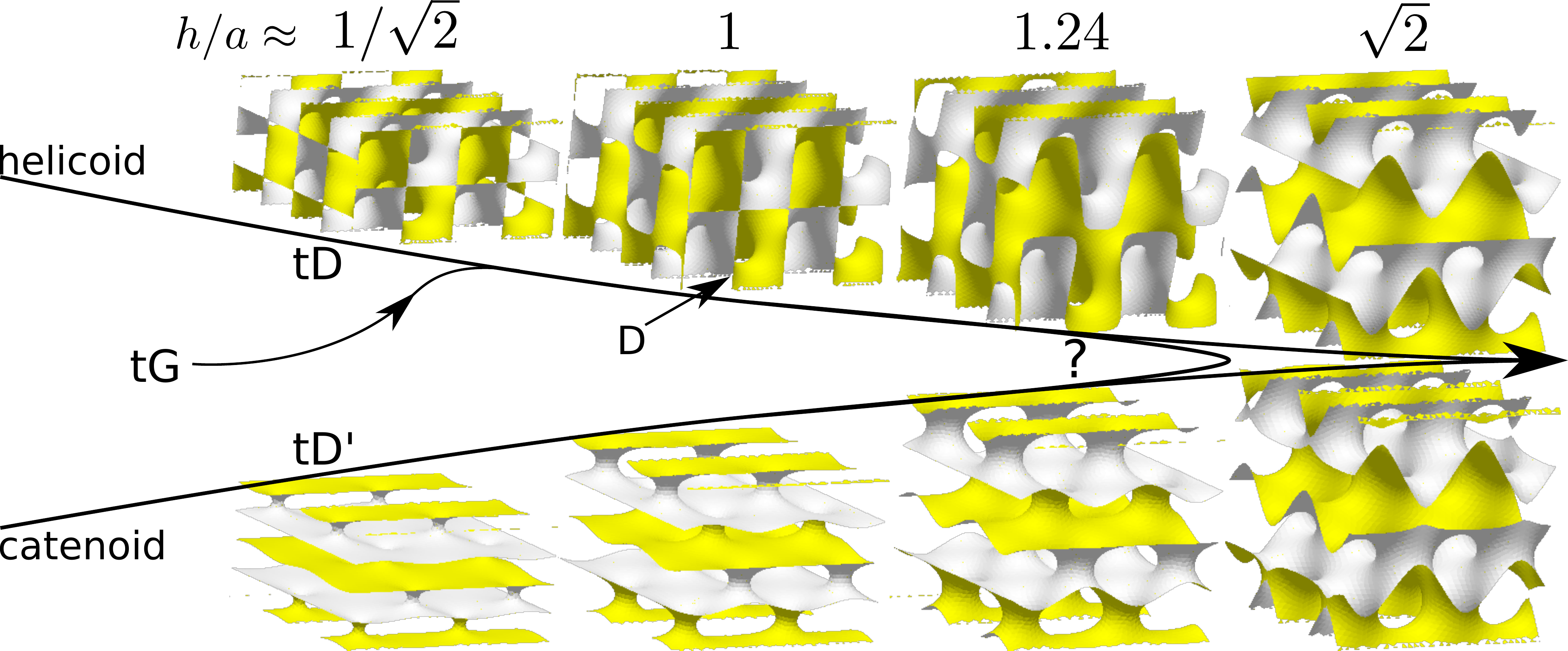}
	\caption{
		\label{fig:tD} Surface Evolver simulations of the catenoid branch ($\tD'$)
		and the helicoid branch ($\tD$) of the extended $\tD$ family, with
		different ratio of the vertical ($h$) and horizontal ($a$) lattice
		parameters.  Surfaces in this figure has been obtained by minimizing the
		Willmore energy to practically $0$.  Extreme caution has been taken for the
		catenoid branch ($\tD'$) to avoid crash.  The $\D$ surface is highlighted.
		Position where $\tG$ family joins $\tD$ is also sketched.  Question mark at
		the junction indicates our ignorance.
	}
\end{figure}

By extremely careful manipulations, we managed to produce both families in
Surface Evolver; see Figure~\ref{fig:tD}.  The upper half of the figure is the
$\tD$ family.  We use $\tD'$ to denote the Traizet family, shown in the lower
half of the figure.  The $\tD'$ family is produced in Surface Evolver by
explicitly connecting small necks between horizontal planes.  Extreme caution
has been taken to avoid crash.

The two families can be understood in the following way.  In general, if two
horizontal planes are close to each other, a catenoid neck between them can
have two possible sizes.  This phenomenon is famous for the catenoid, and also
mentioned when we construct $\rPD$ surfaces.  Our experiments shows that
smaller necks correspond to the $\tD'$ family, and larger necks correspond to
the $\tD$ family.

As we increase the distance between $(100)$-planes, the $\tD$ and $\tD'$
families seem to merge into one.  We do not fully understand the behavior of
surfaces near the junction.  We propose the term \emph{extended $\tD$ family}
for the union of the two families.

\medskip

We end this paper with some discussions on the twinning of $\mG$ surface.  Near
Traizet's catenoid limit, the $\mG$ surfaces can certainly be twinned about
horizontal planes; the force balancing condition off the boundary is covered
by~\eqref{eq:mG}, and the condition on the boundary is trivial.  But in view of
the observed $\G$ twin, we are more interested in twinning $\mG$ surfaces about
vertical planes parallel to $T_1$.  We performed numerical calculations with
$\delta=1,2,3$.  However, if we insist that the twin boundaries pass through
the catenoid necks, the only balanced configurations turn out to
be~\eqref{eq:mG} with $T_1=1$, $T_2=bi$ or $T_2=1+bi$, where $b=\sqrt{8/9}$.
This is not exactly what we observed in the $\G$ twins; thus our understanding
of $\G$ twins, and more generally $\mG$ twins, is still very limited.  Note
that the case $T_2=bi$ is an orthorhombic deformation of the $\D$ surface (oD)
as we have observed, which should not be a surprise any more in view of
previous discussions on the $\D$ surface.

For a concrete example of calculation, the configuration for $\delta=1$ is (up to
scaling)
\begin{gather*} \label{eq:oD}
	\begin{split}
		&N=4,\quad m_k=1\quad (k=0\dots3),\\
		&p_{0,1}=0,\quad p_{1,1}=1/2,\\
		&p_{2,1}=x+1/2+bi/2,\quad p_{3,1}=x+bi/2,\\
		&T_1=1,\quad T_2=bi,\quad T_3=0.
	\end{split}
\end{gather*}
The only solutions for $F_{0,1}=0$ are $x=0$ and $x=1/2$.


\end{document}